\documentclass[12pt]{article}
\title{{\bf Rational Hypergeometric Functions}}
\author{Eduardo Cattani, Alicia Dickenstein and Bernd Sturmfels}
\date{November 4, 1999}

\usepackage{amsthm}
\usepackage{amsmath}
\usepackage{amsfonts} 

\newcommand{\baseRing}[1]{\ensuremath{\mathbb{#1}}}
\newcommand{\Z}{\baseRing{Z}}
\newcommand{\R}{\baseRing{R}}
\newcommand{\C}{\baseRing{C}}
\newcommand{\N}{\baseRing{N}}
\newcommand{\Q}{\baseRing{Q}}
\newcommand{\CP}{\baseRing{P}}

\def\pd#1{ \partial_{#1} }

\theoremstyle{plain}
\newtheorem{theorem}{Theorem}[section]
\newtheorem{lemma}[theorem]{Lemma}
\newtheorem{corollary}[theorem]{Corollary}
\newtheorem{prop}[theorem]{Proposition}
\newtheorem{conjecture}[theorem]{Conjecture}

\theoremstyle{definition}
\newtheorem{definition}[theorem]{Definition}
\newtheorem{remark}[theorem]{Remark}

\newtheorem{example}[theorem]{Example}

\numberwithin{equation}{section}

\newcommand{\Script}[1]{\ensuremath{{\cal{#1}}}}

\newcommand{\RR}{\Script{R}}

\def\D{D}

\begin{document}
\maketitle
\begin{abstract}
Multivariate hypergeometric functions associated with toric
varieties were introduced by Gel'fand, Kapranov and Zelevinsky.
Singularities of such functions are discriminants, that is, 
divisors projectively dual to torus orbit closures.
We show that most of these potential denominators never appear in rational 
hypergeometric functions. We conjecture that the denominator of any rational
hypergeometric function is a product of resultants,  that is, a product of 
special discriminants arising from Cayley configurations.
This conjecture is proved for toric hypersurfaces and for toric varieties of 
dimension at most three. 
Toric residues are applied to show that
every toric resultant appears in the denominator of some rational
hypergeometric function.

\end{abstract}

\setcounter{section}{0}
\section{Introduction}
\label{section:1}
\setcounter{equation}{0}

Which rational functions  in $n$ variables
are hypergeometric functions? \linebreak
Which denominators 
appear in such rational hypergeometric functions?
Our aim is to answer these questions for the 
multivariate hypergeometric functions
introduced by Gel'fand, Kapranov and 
Zelevinsky \cite{gkz89,gkz90,sst2}.
These functions are defined by a system
of linear partial differential equations,  associated to
any integer   $d \times s$-matrix $A  = (a_{ij})$
and any complex vector $\beta \in \C^d$:

\begin{definition}\label{def:hypergeom}
The {\it $A$-hypergeometric system of degree  $\beta\in\C^{d }$}
 is  the left ideal 
$H_A(\beta)$ in the Weyl algebra
$ \C\langle x_1,\dots,x_s,\pd 1,\dots,\pd s\rangle$
 generated by
\begin{eqnarray}
\label{higherorder}
\hbox{the toric operators} \quad
\partial^u - \partial^v \quad \hbox{for}\quad 
 u,v\in \N^s    \quad \hbox{with}\quad A\cdot u=A\cdot v, 
\phantom{wow} \\
\label{firstorder}
\hbox{and the Euler operators} \quad
\sum_{j=1}^s a_{ij} x_j \partial_j - \beta_i
\quad \hbox{for} \quad i=1,\dots,d. \phantom{wowowow}
\end{eqnarray}
A function $f(x_1,\dots,x_s)$, holomorphic on an open
set $U\subset \C^s$, is said to be $A$-hypergeometric of
degree $\beta$ if it is annihilated by the left ideal $H_A(\beta)$.
\end{definition}

Throughout this paper we use the multi-exponent notation
$\partial^u  = \prod_{j=1}^s \partial_j^{u_j} \!$.
We shall assume that the
rank of the matrix $A$ equals $d$, 
the column vectors $a_j$ of $A$ are distinct, and
the vector $(1,1,\dots,1)$ lies in the row span of $A$.

These hypotheses greatly simplify our exposition, but our
main results remain valid without them.
The last hypothesis means that the toric ideal
$$
I_{A}\  :=\  \langle\, \xi^u - \xi^v \,\, : \,\,
A\cdot u = A\cdot v \,\rangle 
\quad \subset \quad \C[\xi_1,\dots,\xi_s]
$$
is  homogeneous with respect to total degree
and defines a projective toric variety  $X_A\subset \CP^{s-1}$, and
the columns of $A$ represent a configuration $\{a_1,\dots,a_s\}$
of $s$ distinct points in affine $(d-1)$-space. This condition
ensures that the system $H_A(\beta)$ has only regular singularities
(\cite{gkz89}, \cite[Theorem~2.4.11]{sst2}). A detailed analysis of
the non-regular case was carried out by Adolphson
\cite{adolphson}.

The system $H_A(\beta)$ is always holonomic. Its holonomic rank $r_A(\beta)$
coincides with the dimension of the space
of local holomorphic solutions in  $\C^s \backslash
{\rm Sing}(H_A(\beta))$.
If $I_A$ is Cohen-Macaulay or $\beta$ is generic in $\C^d$, then
\begin{equation}\label{rankequation}
r_A(\beta) \quad  = \quad  {\rm degree}(X_A)
\,\,\, = \,\,\,  {\rm vol}( {\rm conv}(A)),
\end{equation}
the normalized volume of  the lattice polytope
${\rm conv}(A) = {\rm conv} \{a_1,\dots,a_s\}$.
The inequality $\, r_A(\beta)   \geq  {\rm vol}( {\rm conv}(A)) \,$
holds without any assumptions on $A$ and $\beta$.
See \cite{adolphson}, \cite{gkz89}, \cite{sst2}  for proofs and details.
If $d=2$, i.e.~when $X_A$ is a curve, then (\ref{rankequation}) holds for
all $\beta\in\C^2$ if and only if $I_A$ is Cohen-Macaulay \cite{cdd}.

The irreducible components of ${\rm Sing}(H_A(\beta))$ are the
hypersurfaces defined by  the $A'$-discriminants $D_{A'}$,
where $A'$ runs over {\it facial subsets} of $A$, or, equivalently,
$X_{A'}$ runs over closures of torus orbits on $X_A$.
The {\it $A$-discriminant} $\D_{A}$ is the
irreducible polynomial defining the dual variety of the toric
variety $X_{A}$, with the convention $\D_{A} = 1$ if that
dual variety is not a hypersurface; see
\cite{adolphson}, \cite{gkz89}, \cite{gkzbook}.
Note that for a singleton
$\,A' = \{a_j\}\,$ we have
$\,\D_{A'} = x_j$.

Consider any rational  $A$-hypergeometric function of degree $\beta$,
\begin{eqnarray}
\label{rationalfct}
f(x_1,\dots,x_s) \quad = \quad
\frac{ P(x_1,\dots,x_s)}{ Q(x_1,\dots,x_s) }, 
\end{eqnarray}
where  $P$ and $Q$ are relatively prime polynomials. 
The denominator equals
\begin{equation}
\label{denomproduct}
 Q(x_1,\dots,x_s) \quad \quad = \quad
\prod_{A'} 
\D_{A'}(x_1,\dots,x_s)^{i_{A'}} ,
\end{equation}
where  $A'$ runs over facial subsets of $A$
and the $\,i_{A'} \,$ are non-negative integers.

Our long-term goal is to classify all
rational  $A$-hypergeometric functions.
For toric curves this was done in \cite{cdd}:
if $d=2$, every rational $A$-hypergeometric function
is a Laurent polynomial. Here we
generalize this result to higher-dimensional toric varieties,
by studying rational $A$-hypergeometric functions
which are not Laurent polynomials. We note that, by 
\cite[Lemma~3.4.10]{sst2}, $A$-hypergeometric
polynomials exist for all toric varieties $X_A$ .

We call the matrix $A$ {\it gkz-rational} if the
$A$-discriminant $\D_A$ is not a monomial
and appears in the denominator  (\ref{denomproduct}) of some
rational $A$-hypergeometric function (\ref{rationalfct}).  
The smallest example  of a gkz-rational configuration is
\begin{equation}
\label{GaussSquare}
A \quad = \quad \Delta_1 \times \Delta_1 \quad = \quad
  \left(
\begin{array}{cccc}
1 & 0 & 1 & 0 \\
0 & 1 & 0 & 1 \\
0 & 1 & 1 & 0 
\end{array}
\right),
\end{equation}
since $\, 1 / (x_1 x_2 - x_3 x_4) \,$ is $A$-hypergeometric
of degree $\, \beta = (-1,-1,-1)^T$.
Note that (\ref{GaussSquare})
encodes the {\sl Gauss hypergeometric function} ${}_2 F_1 \,$
 \cite[\S 1.3]{sst2}.

More generally, the product of simplices
$\,A = \Delta_p \times \Delta_q$ 
is gkz-rational if and only if $p=q$.
The Segre variety $X_A = \CP^q \times \CP^q$ is projectively dual to the
$(q \! + \!1) \times (q \! + \!1)$-determinant, and the
reciprocal of this determinant is a rational $A$-hypergeometric function.
Consider by contrast the configuration  $\, A = 2 \cdot \Delta_q $.
The toric variety $X_A$ is the
quadratic Veronese  embedding of $\CP^q$, whose
projectively dual hypersurface is
 the discriminant of a quadratic form,
\begin{equation}
\label{DiscrQuadForm}
\D_A \quad = \quad
\det
  \left(
\begin{array}{ccccc}
2 x_{00} & x_{01} & x_{02} & \cdots & x_{0q} \\
x_{01}   & 2 x_{11} & x_{12} & \cdots & x_{1q} \\
x_{02}   &  x_{12}  & 2 x_{22} & \cdots & x_{2q} \\
 \vdots &  \vdots &  \vdots & \ddots & \vdots \\
x_{0q}   & x_{1q} & x_{2q} & \cdots & 2 x_{qq} 
\end{array}
\right).
\end{equation}
Theorem \ref{MainThm} below implies that 
the classical (``dense'') discriminants
such as (\ref{DiscrQuadForm}) do not appear in the
denominators of rational hypergeometric functions.
In other words, multiples of simplices, $A = r \cdot \Delta_q$,
are  never gkz-rational. 

In Section 2 we resolve the case of {\it circuits}, that is, matrices
$A$ whose kernel is spanned by a single
vector $b = b_+ - b_- \in \Z^s$. We call $A$ 
{\it balanced} if the positive part $b_+$ is a coordinate
permutation of the negative part $b_-$, and we show
that $A$ is balanced if and only if $A$ is gkz-rational.
In Section 3 we study arbitrary configurations $A$, and
we prove the following theorem.

\begin{theorem}\label{MainThm}
If the configuration $A$ contains an unbalanced circuit 
which does not lie in any proper facial subset of $A$,
then $A$ is not gkz-rational.
\end{theorem}

This implies that  gkz-rational configurations are rare; for instance,
they have no interior points. Hence, reflexive polytopes \cite{batyrev}
are not gkz-rational, and
sufficiently ample embeddings of any toric variety are not gkz-rational.

In order to formulate a
conjectural characterization of gkz-rational configurations, we recall
the following construction from \cite{gkzbook}.
 Let $A_0,A_1,\dots,A_r$ be vector configurations in $\Z^r$.
Their {\it Cayley configuration} is defined as
\begin{equation}
\label{Cayley}
  A \,\,\, = \,\, \,
 \{e_0\} \! \times A_0 \! \, \, \cup \,\, 
 \{e_1\} \! \times A_1 \! \, \, \cup \,\, 
\cdots   \, \, \cup \,\,
\{e_{r}\}   \! \times \! A_{r}
\,\, \subset \,\, \Z^{r+1} \times \Z^r , 
\end{equation}
where $e_0,\dots, e_r$ is the standard basis of $\Z^{r+1}$.

We call $A$ {\it essential} if the Minkowski sum
$\,\sum_{i \in I} A_i \,$ has affine dimension at least
$|I|$ for every proper subset $I$ of $\{0,\dots,r\}$.
Cayley configurations are very special.
 For instance, a configuration $A$ in the plane
 $(d=3)$ is a Cayley configuration if and only if 
$A$ lies on two parallel lines;
such an $A$ is  essential
if and only if each  line contains at least two points.

\begin{conjecture} \label{conj4}
An arbitrary configuration $A$ is gkz-rational if and only if $A$ is 
affinely isomorphic to an essential Cayley configuration (\ref{Cayley}).
\end{conjecture}

This conjecture can be reformulated as follows.
The discriminant  of an essential Cayley configuration coincides with the 
{ sparse resultant} $R_{A_0,A_1,\dots,A_r}$; see \cite[\S 8.1.1]{gkzbook}.
That resultant  characterizes the solvability of a
sparse polynomial system $\, f_0 = f_1 =  \cdots = f_r = 0 \,$
with support $(A_0,A_1,\dots,A_r)$, 
$$ f_i(t_1,\dots,t_r) \quad = \quad \sum_{a \in A_i} x_a  t^a ,
\qquad i  = 0,1,\dots, r.      $$
By Corollary \ref{5.2}, Conjecture \ref{conj4}  is
equivalent to the following:

\begin{conjecture} \label{conj4prime}
A discriminant $D_A$ appears in the denominator of
a rational $A$-hypergeometric function if and only if $D_A$
is a resultant $\,R_{A_0,A_1,\dots,A_r}$.
\end{conjecture}

Being a resultant among  discriminants is being a needle in a haystack.
None of the univariate or classical discriminants 
such as (\ref{DiscrQuadForm}) are resultants.
On the other hand, consider
two triples of equidistant points on parallel lines,
\begin{equation}
\label{Scroll}
 A \quad = \quad
\{e_0\}  \! \times \! A_0  \, \, \cup \,\, 
\{e_1\}  \! \times \! A_1  
\quad = \quad
\left(
\begin{array}{rrrrrr}
1 & 1 & 1 & 0 & 0 & 0 \\
0 & 0 & 0 & 1 & 1 & 1 \\
0 & 1 & 2 & 0 & 1 & 2 
\end{array}
\right).
\end{equation}
This is the Cayley configuration of
$  A_0 = A_1 = \{ 0 , 1 , 2 \}$.
The variety $X_A$ is a rational normal scroll in $\CP^5$.
Its discriminant $ D_A$ is the Sylvester resultant
$$
R_{A_0,A_1} \,  = \,
x_1^2 x_6^2-x_1 x_2 x_6 x_5-2 x_1 x_3 x_4 x_6+x_1 x_3 x_5^2
+x_2^2 x_4 x_6-x_2 x_3 x_4 x_5+x_3^2 
x_4^2 $$
of the quadrics
$ \,F_0  = x_1 u_1^2+ x_2 u_1 u_2 +x_3 u_2^2 \,$ and
$\,F_1 =   x_4 u_1^2 + x_5 u_1 u_2 + x_6 u_2^2 $.
The following theorem is the second main result in this paper.

\begin{theorem}
\label{secondmainthm}
The if directions of Conjectures  \ref{conj4} and \ref{conj4prime} 
hold. The only-if directions hold for
$\, d \leq 4 $, that is,  for toric varieties $X_A$
of dimension $\leq 3$.
\end{theorem}

The proof of the only-if direction is given in Section 4.
It consists of a detailed combinatorial case analysis based on
Theorem \ref{MainThm}. The proof of the if direction, given in
Section 5, is based on the  notion of {\it toric residues} introduced by 
Cox \cite{cox},  and on our earlier results 
in \cite{tokyo} about their denominators.

An example of a toric residue is the rational  $A$-hypergeometric function 
\begin{equation}
\label{anintegral}
 \frac {1}{(2\pi i)^2} \int_\Gamma \ 
\frac {u_1 \,u_2 }{F_0(u_1,u_2) \cdot F_1(u_1,u_2)}\,du_1\wedge du_2 
\quad = \quad \frac{x_1 x_6 - x_3 x_4}{ R_{A_0,A_1}}.
\end{equation}
Here $A$ is the configuration
(\ref{Scroll}) and
 $\Gamma$ is a suitable $2$-cycle  in $ \C^2$.
Such integrals can be evaluated 
by a single Gr\"obner basis computation; see \cite{global}.

\bigskip

\section{Circuits}

\label{section:2}
\setcounter{equation}{0}

We fix a configuration $A$ which is
a circuit, that is, $A$ is a $d \times (d+1)$-matrix whose integer
kernel is spanned by a vector $\,b = (b_0,b_1,\dots,b_d)\,$
all of whose coordinates $b_i$ are non-zero.  After relabeling, 
we may assume
\begin{equation}
b_j > 0 \  \hbox{for}\ j=0,\dots,m-1\  \quad \hbox{and} \quad
b_j < 0 \  \hbox{for}\ j=m,\dots,d,
\end{equation}
so that $b_+=(b_0,\dots,b_{m-1},0,\dots,0)$ and 
$b_-=(0,\dots,0,-b_{m},\dots,-b_d)$.
The toric variety $X_A$ is a hypersurface
in $\CP^d$, defined by the principal ideal
$$ I_A \quad = \quad \langle \, 
\xi^{b_+} - \xi^{b_-} \,\rangle
 \quad = \quad \langle\,
\xi_0^{b_0} \cdots \xi_{m-1}^{b_{m-1}}  \, - \,
\xi_{m}^{-b_m} \cdots \xi_d^{-b_d} \, \rangle .  $$
In this section we shall prove our main conjecture for the 
case of circuits.

\begin{theorem} \label{Conj4TrueforCircuits}
Conjectures \ref{conj4} and
\ref{conj4prime} are true for toric hypersurfaces.
\end{theorem}

A function $f(x_0,x_1,\dots,x_d)$ is $A$-hypergeometric
if it is $A$-homogeneous  (satisfies (\ref{firstorder}) for some $\beta$)
and annihilated by the
homogeneous toric operator
\begin{equation}
\label{Boperator}
\partial^{b_+} - \partial^{b_-} \quad = \quad
\partial_0^{b_0} \cdots \partial_{m-1}^{b_{m-1}}  \, - \,
\partial_{m}^{-b_m} \cdots \partial_d^{-b_d}  . 
\end{equation}
The order $\rho$ of this operator equals the holonomic rank of $H_A(\beta)$:
$$
\rho := \ = \
b_0 + \cdots + b_{m-1} \  = \ 
- b_m - \cdots - b_d  \  = \   {\rm vol}(
{\rm conv}(A))\  = \ r_A(\beta).
$$
This holds for all $\beta \in \C^d$
since the principal ideal
$I_A = \langle \xi^{b_+} - \xi^{b_-} \rangle$ is Cohen-Macaulay.
The toric hypersurface $X_A$ is 
projectively self-dual. Indeed, by \cite[Proposition 9.1.8]{gkzbook},
the $A$-discriminant of the circuit $A$ equals
\begin{equation}
\label{circuitdisc}
 D_A \quad = \quad  x^{b_-} - \lambda x^{b_+} \qquad
\hbox{ where $\,\lambda = (-1)^\rho  b_-^{b_-}/ b_+^{b_+}$.}
\end{equation}
Recall that the circuit $A$ is {\it balanced} if
$d = 2 m-1$ and, after reordering if 
necessary, $b_i = -b_{m+i}$ for $i=0,\dots,m-1$. Otherwise, we 
call $A$ {\it unbalanced}.
Note that the configuration (\ref{GaussSquare}) is a balanced circuit 
with $b = (1,1,-1,-1)$.

\begin{lemma} Let $A$ be a balanced circuit. Then 
the rational function  $1/ D_A$ is $A$-hypergeometric.
\end{lemma}

\begin{proof} Balanced implies 
$\, b_-^{b_-} = b_+^{b_+}\,$ and  $\lambda = (-1)^\rho$.
Consider the expansion
$$ \frac{1}{D_A} \quad = \quad
x^{-b_-}  \cdot \frac{1 }{ 1 - (-1)^\rho x^b}
\quad = \quad \sum_{n = 0}^\infty (-1)^{\rho n} x^{n b_+ - (n+1) b_-} .$$
For this series to be
annihilated by  (\ref{Boperator}) it is
necessary and sufficient that
$$
\prod_{i=0}^{m-1}  \prod_{j=1}^{b_i} (n b_i+j) \,\,\, = \,\,\,
\prod_{i=m}^s  \prod_{j=1}^{-b_i} (n (-b_i)+j) 
\quad \hbox{for all $n \geq 0$.}
$$
This identity holds if and only if the circuit $A$ is balanced.
\end{proof}

This lemma implies that balanced circuits are gkz-rational.
The main result in this section is the following converse to this
statement. For an arbitrary configuration $A$, we say that $A$ 
is {\it weakly
gkz-rational} if there exists a rational $A$-hypergeometric 
function which is not a Laurent polynomial.

\begin{theorem}\label{th:2.1}
Let $A$ be a circuit in $\Z^d$. Then the following are equivalent:
(1) $A$ is balanced; \ 
(2)  $A$ is gkz-rational; \ 
(3) $A$ is weakly gkz-rational.
\end{theorem}

\begin{proof} The implication from (1) to (2) follows from the
previous lemma. The equivalence  of (2) and (3) holds because
every proper facial subset $A'$ of  $A$ is affinely independent.
Hence the only non-constant $A'$-discriminants arising from facial
subsets $A'$ arise from vertices $A' = \{a_j\}$, in which case
$D_{A'} = x_j$.

It remains to prove the implication from (2) to (1). Suppose
that $A$ is gkz-rational.
Consider a non-Laurent rational $A$-hypergeometric function
and expand it as a Laurent series with respect to increasing
powers of $x^b$.  It follows from the
results in \cite[\S 3.4]{sst2} that this series is the sum
of a Laurent polynomial and a 
{\it canonical $A$-hypergeometric series} of the following form:
\begin{equation}\label{eq:2.17}
F (x) \quad = \quad x^c \cdot \sum_{n=0}^\infty \,
 (-1)^{\rho n} \ 
\frac
{\prod_{j \geq m}(-c_j-nb_j-1)!}{\prod_{j < m}(c_j + n b_j)!} \ x^{nb}.
\end{equation}
Here $c = (c_0,c_1,\dots,c_d)$ is a suitable integer vector. 
The series $F(x)$ represents a rational function.
We may view the series on the right-hand side of  (\ref{eq:2.17})
as defining a rational function of 
the single variable $\,t = x^b  ={x^{ b_+}}/{x^{ b_-}}$:
\begin{equation}\label{eq:2.8}
\varphi(t)
\quad = \quad \sum_{n=0}^\infty \,
 (-1)^{\rho n} \ 
\frac
{\prod_{j\geq m}(-c_j-nb_j-1)!}{\prod_{j < m}(c_j + n b_j)!} \,
{t^n}
\end{equation}

The $A$-discriminant equals
$D_A   =  x^{b_-}(1- \lambda t)$
where $\lambda = (-1)^\rho \, b_-^{b_-}/ b_+^{b_+}$.  This implies that
the rational function $\varphi(t)$ may be written as a quotient
$$\varphi(t) \ =\  \frac {P(t)}{(1- \lambda t)^{k+1}},$$
where $P (t)$ is a polynomial and $k \in \N$.
It follows from \cite[Corollary 4.3.1]{stanley} that the 
coefficients of the series (\ref{eq:2.8}) must be of the form
$\lambda^n$ times a polynomial in $n$. That is,
the following expression is a polynomial in $n$:
$$\gamma(n) \ :=\  \lambda^{-n}\cdot 
\frac
{\prod_{j \geq m}(-c_j-nb_j-1)!}{\prod_{j <  m}(c_j + n b_j)!}$$ 
The rational function
$\,\mu(z) := \gamma(z+1)/\gamma(z)\,$ satisfies
the following general identity  \cite[Lemma 2.1]{singer}
for any fixed complex number $z_0$:
\begin{equation}\label{eq:2.10}
\sum_{\alpha \in z_0 + \Z}  \ {\rm ord}_\alpha(\mu) \,\,\, = \,\,\, 0.
\end{equation}
Our rational function
$\mu(z)$  has its poles among the points
$$-(\frac {c_j}{b_j} + \frac {1}{b_j}), \,
-(\frac {c_j}{b_j} + \frac {2}{b_j}),\, \dots,\,
-(\frac {c_j}{b_j} + 1)\ ;\ \ j=0,\dots,m-1,$$
and its zeroes among 
$$-\frac {c_j}{b_j},\,
-(\frac {c_j}{b_j} + \frac {1}{b_j}),\, \dots,\,
-(\frac {c_j}{b_j} + \frac {-b_j-1}{b_j})\ ;\ \ j=m,\dots,d.$$
We may assume 
$b_0 = {\rm max}\{b_j\ ; j =0,\dots,m-1\}$ and 
$-b_m = {\rm max}\{-b_j\ ; j =m,\dots,d\}$.  Suppose now
that $b_0 > -b_m$.  Then, $\mu(z)$ has a
pole at a point $p/b_0$ with $p$ and $b_0$ coprime, but 
since none of the zeroes may be of this form, this
contradicts (\ref{eq:2.10}).  A symmetric argument
leads to a contradiction if we assume $b_0 < -b_m$.
This implies that $b_0 = -b_m$ and therefore
$$\gamma(n) 
\cdot \frac {(c_0 + n b_0)!}{(-c_m-n b_m-1)!} $$
is also rational function of $n$.  Consequently, we can iterate our
argument to conclude that, after reordering,
$b_i = -b_{m+i}$ for all $i=0,\dots,m-1$.
\end{proof}

\begin{remark}
The above results imply that a circuit $A$ is gkz-rational
if and only if the specific rational function $1/ D_A$ is
A-hypergeometric. The same statement is false for non-circuits.
For instance, for the gkz-rational configuration in (\ref{Scroll}),
the function  $1/D_A = 1/R_{A_0,A_1}$ is not $A$-hypergeometric.
\end{remark}

Let us now return to the result stated at the beginning of
this section.

\begin{proof}[Proof of Theorem \ref{Conj4TrueforCircuits}]
Theorem \ref{th:2.1} and the lemma below imply
Conjecture \ref{conj4}.
The equivalence of 
Conjectures \ref{conj4} and
\ref{conj4prime} will be shown in Section 5.
\end{proof}

\begin{lemma}\label{balcirc} A circuit $A$ is balanced if and only if
it is affinely isomorphic to an essential Cayley configuration (\ref{Cayley}).
\end{lemma}

\begin{proof}
We first prove the if-direction. Let $A$ be an essential
Cayley configuration which is a circuit.
Each  $A_i$ must consist of a pair of vectors in $\Z^r$,
so that $A$ becomes an $(2r+1) \times (2r+2)$-matrix. 
The first $r+1$ rows of $A$ show that the kernel of $A$ is
spanned by a vector $\, b = (b_0,-b_0,b_1,-b_1,\dots,b_r,-b_r)$.
This means that $A$ is balanced.
Conversely, if $A$ is balanced then we can apply
left multiplication by an element of $GL(d,\Q)$ to get isomorphically
$$ A \quad = \quad \left( \begin{array}{cc}
 I_m & I_m \\ 0& \tilde A
\end{array} \right) $$
where $\tilde A$ is an $(m-1)\times m$ 
integral matrix of rank $m-1$. 
By permuting columns we see that $A$ is an essential
Cayley configuration for $m=r+1$.
\end{proof}

\section{The General Case} 

In this section we prove  Theorem~\ref{MainThm}.
A configuration $A$ is called
{\em non-degenerate\/} if
the $A$-discriminant $D_A$ is neither equal to $1$
nor a variable. Circuits are
non-degenerate by (\ref{circuitdisc}).
Recall that $D_A$ is a variable if and only if
$A$ is a point.
A subconfiguration $B \subseteq A$
is called {\it spanning} if $B$ is not contained
in any proper facial subset of $A$.   
If the dimension of
 $B$ is equal to the dimension of $A$
then $B$ is spanning, but the converse is not true.
For instance, the vertex set of an octahedron
contains spanning circuits but no full-dimensional circuits.

The condition $D_A = 1$ 
means that  the dual variety to
the toric variety $X_A$ is not a hypersurface. 
No combinatorial characterization of this condition is presently known.
A necessary condition is given in the following proposition.
That condition is not sufficient: the skew prisms in (\ref{SixPoints})
contain no spanning circuit but $D_A \not= 1$. Note that
$D_A=1$ for the regular prism $A= \Delta_1 \times \Delta_2$.

\begin{prop}\label{nondeg}
If  $A$ contains a subconfiguration $B$ which is
spanning  and non-degenerate then $A$ is non-degenerate.   
In particular, $A$ is non-degenerate
if it contains a spanning circuit.
\end{prop}

\begin{proof}
Proceeding by induction, it suffices to consider the case
when $B$ is obtained from $A$ by removing a single point, say,
$B = A \backslash \{a_s\}$.  
Since $B$ is not contained in any face of $A$, 
and $B$ is a facial subset of itself,
the following lemma tells us that the $B$-discriminant $D_B$
divides $ \,D_{A}|_{x_s=0}$. Since
$D_B$ is not a monomial, this implies that  $D_A$ is not a monomial.
\end{proof}

\begin{lemma}\label{setzero}
Let $a_s \in A$, $x_s $ the corresponding variable, and
$B' $ a facial subset of $B = A \backslash \{a_s \}$
which does not lie in any proper facial subset of $A$.
Then the $B'$-discriminant $D_{B'}$ divides the
specialized $A$-discriminant $ \,D_{A}|_{x_s=0}$.
\end{lemma}

\begin{proof}
Let $f = \sum_{j=1}^s x_j t^{a_j} $ be a generic polynomial
with support $A$. 
By \cite[Theorem~5.10]{ksz},
the {\it principal $A$-determinant}
 is the specialization
$$ E_A \quad = \quad
R_A \bigl(
 t_1 {\frac {\partial f} {\partial t_1}},\dots,
 t_d {\frac {\partial f} {\partial t_d}} \bigr) , $$
where $R_A$ denotes the {\it $A$-resultant}; 
see \cite[\S 8.1]{gkzbook}.
The irreducible factorization of
the principal $A$-determinant ranges
over the facial subsets $A'$ of $A$,
\begin{equation}
\label{PrincipalDet}
 E_A \quad = \quad
\prod_{A'} D_{A'}^{m_{A'}},
\end{equation}
where $m_{A'}$ are certain positive integers
\cite[Theorem 10.1.2]{gkzbook}.

Let $w \in \Z^s$ be the weight vector with
$w_s = -1$ and $w_j = 0$ for $j \not= s$.
The initial form of the  principal
$A$-determinant with respect to $w$ can
be factored in two different ways:
$$
 in_{w}(E_A) \quad
= \quad   
\prod_{A' } (in_{w} D_{A'})^{m_{A'}} 
= \quad   
\prod_{C \,{\rm facet} \, {\rm of} \, \Delta_w }  (E_C)^{n_{C}} .
$$
Here $\Delta_w$ is the coherent polyhedral subdivision of $A$
defined by $w$ and the $n_C$ are certain positive integers.
The first formula comes from  (\ref{PrincipalDet}) and the second
formula comes from \cite[Theorem 10.1.12]{gkzbook}.
Since $B'$ is a facial subset of $A \backslash \{a_s\}$, 
it is also a cell of the subdivision $\Delta_w$, and hence
$D_{B'}$ divides $E_C$ for the facet $C= A\backslash  \{a_s\}$ of $\Delta_w$.
We conclude that $D_{B'}$ divides
$\, in_w(D_{A'})\,$ for some facial subset $A'$ of $A$.
If $D_{B'} \not= 1$, this implies that $B' \subseteq A'$ 
because $D_{B'}$ involves all the variables associated with
points in $B'$.
By our hypothesis, the only facial subset of $A$ which contains
$B'$ is $A$ itself. Therefore $D_{B'}$ divides
$\, in_w(D_A)\, = \, D_A|_{x_s = 0}$.
\end{proof}

We also need the following lemma from commutative algebra
whose proof was shown to us by Mircea Musta\c{t}\v{a}:

\begin{lemma}\label{mircea}
Let ${\cal R}$ be a unique factorization domain with field of fractions $K$, 
and  let $\,f (t) = \sum_{i=0}^m a_i \cdot t^i \,$ and $\,g (t) = 
\sum_{j=0}^n b_j \cdot t^j \,$ be  relatively prime elements in the 
polynomial ring ${\cal R}[t]$. Assume that $b_0 \not=0 $ and consider
the Taylor series expansion of the ratio $f/g$:
$$  {\frac {f (t)}{g(t)}} \  = \ 
\sum_{\ell=0}^\infty c_\ell \cdot t^\ell \quad \hbox{ in $\, K[[t]]$}. $$
If all the Taylor coefficients 
$\,c_\ell \,$ lie in $ {\cal R}$, then $b_0 $ is a unit in ${\cal R}$.
\end{lemma}

\begin{proof}
Let $p$ be any prime element in ${\cal R}$. We must show that
$p$ does not divide $b_0$.
Consider the localization $\,{\cal R}[t]_{\langle p, t \rangle}\,$
of ${\cal R}[t]$ at the prime ideal $\langle p, t \rangle $.
The power series ring ${\cal R}[[t]]$ is the completion 
of the local ring  $\,{\cal R}[t]_{\langle p, t \rangle}\,$
with respect to the $\langle t \rangle$-adic topology.
By assumption, the polynomial $f$ lies in the principal ideal generated
by $g$ in ${\cal R}[[t]]$. The basic flatness property of completions, 
as stated in \cite[\S 8, p.63]{matsumura},
implies that $f$ lies in the principal ideal generated
by $g$ in $\,{\cal R}[t]_{\langle p, t \rangle}$.
Since $f$ and $g$ are relatively prime in ${\cal R}[t]$,
we conclude that $g$ is a unit in 
$\,{\cal R}[t]_{\langle p, t \rangle}\,$
and so $b_0$ is not divisible by $p$. 
\end{proof}

We are now prepared to prove the first theorem stated in the introduction.

\begin{proof}[Proof of Theorem \ref{MainThm}]
Suppose $A = \{a_1,\dots,a_s\}$ is a gkz-rational configuration
and let $f = P/Q$ be a rational $A$-hypergeometric function
of degree $\beta\in\Z^d$,
where $P, Q \in \C[x_1,\dots,x_s]$ are relatively prime,
and the $A$-discriminant $D_A$ is not a monomial and divides $Q$.
We claim that any spanning circuit $Z$ of $A$ is balanced.
We shall prove this by induction on $s$. If
$s = d+1$, then we are done by Theorem~\ref{th:2.1}.
We may assume that $A$ is not a circuit and
therefore $Z$ is a proper subset of $A$. Suppose
$a_s \in A \backslash Z$, and set $t=x_s$, 
$\tilde A = \{a_1,\dots,a_{s-1}\}$, $\tilde x = (x_1,\dots,x_{s-1})$.
We may expand the rational $A$-hypergeometric function
$\,f(x) = f(\tilde x;t)\,$ as
\begin{equation}
\label{SeriesForMircea}
f(\tilde x;t) \quad = \quad 
\sum_{\ell\geq \ell_0}\,R_\ell(\tilde x)\cdot t^\ell\,,
\end{equation}
where each $R_\ell(\tilde x)$ is a rational $\tilde A$-hypergeometric
function of degree $\beta - \ell\cdot a_s$.  

Let $A'$ denote the unique smallest facial subset of 
$\,\tilde A =  A \backslash \{a_s\}\,$ which contains 
the circuit $Z$. Then
$Z$ is a spanning circuit in $A'$. 
Proposition~\ref{nondeg} implies that  its
discriminant $D_{A'}$ is not a monomial. 
Lemma~\ref{setzero} implies that $D_{A'}$ divides
$in_t(Q)$, the lowest coefficient of $Q$ with respect to $t$.

We now apply  Lemma~\ref{mircea} to the domain 
${\cal R} = \C[\tilde x,\tilde x^{-1}]_{<D_{A'}>}$, the localization
of the Laurent polynomial ring at the principal prime ideal
$\langle D_{A'} \rangle$. Since $in_t(Q)$ is not a unit in ${\cal R}$,
we conclude that some Taylor coefficient  $R_\ell(\tilde x)$
lies in the field of fractions of ${\cal R}$ but not in ${\cal R}$ itself.
This means that $D_{A'}$ divides the denominator of $R_\ell(\tilde x)$.
We have found a rational $\tilde A$-hypergeometric function
whose denominator contains the non-trivial factor $D_{A'}$.
It follows by induction that the spanning circuit $Z $ of $A'$ is balanced.
\end{proof}

Recall that a configuration $A$ is called 
{\it weakly gkz-rational} if there exists a rational $A$-hypergeometric 
function which is not a Laurent polynomial. It is called
{\it gkz-rational} if the $A$-discriminant $\D_A$ is not a monomial
and appears in the denominator of a rational $A$-hypergeometric function.

\begin{prop}\label{Weakk}
A configuration $A$ is weakly gkz-rational if and only
if some facial subset $A'$ of $A$ is gkz-rational.
\end{prop}

\begin{proof}
If $A'$ is a facial subset of $A$ then every 
$A'$-hypergeometric function $f(x)$ is also $A$-hypergeometric.
Indeed, $f(x)$ is obviously $A$-homogeneous, but it is also
annihilated by the toric operators $\partial^u - \partial^v$ because 
the support of $\partial^u$ lies in $\,\{\partial_i \,:\,a_i \in A'\}\,$
if and only if
the support of $\partial^v$ lies in $\,\{\partial_i \,:\,a_i \in A'\}$.
This proves the if-direction.
For the only-if direction, suppose that $A$ is weakly gkz-rational
and let $f(x) = P(x)/Q(x)$ be a non-Laurent
rational hypergeometric function. There exists a facial subset $A'$
of $A$ such that  $D_{A'}$ is not a monomial and divides $Q(x)$.
Our goal is to show that $A'$ is gkz-rational.
We proceed by induction on the cardinality of  $A \backslash A'$. 
There is nothing to show if $A = A'$. Let $a_s \in A \backslash A'$
and form the series expansion as in (\ref{SeriesForMircea}).
By applying Lemma \ref{mircea} as in the proof of Theorem \ref{MainThm},
we construct a rational $(A \backslash \{a_s\})$-hypergeometric
function whose denominator is a multiple of
the $A'$-discriminant $D_{A'}$.
This proves our claim, by induction.
\end{proof}

We close this section with two corollaries which demonstrate
the scope of Theorem \ref{MainThm}. They show that
gkz-rational configurations $A$ are very special.

\begin{corollary} \label{interior}
A gkz-rational configuration $A$ has no interior point.
\end{corollary}

\begin{proof} 
Let $a_1$ be an interior point of $A$, and let $Z'$ be 
a minimal size subset of $A \backslash \{a_1\}$ which contains
$a_1$ in its relative interior. Then $Z = Z' \cup \{a_1\}$ is a circuit of $A$
which is spanning and not balanced.
\end{proof}

Corollary \ref{interior} can be rephrased,
using Khovanskii's genus formula \cite{askold},
 into the language of algebraic 
geometry as follows. 
If a projective toric variety $X_A$ is gkz-rational, 
then the generic hyperplane section of $X_A$ has arithmetic genus $0$.
Clearly, this fails if $X_A$ is embedded by a sufficiently ample line bundle,
and also in the case of special interest in mirror symmetry 
(see \cite{batyrev}).

\begin{corollary}  The configuration $A$ is not gkz-rational if
$A$ is the set of lattice points in a reflexive polytope, or $A$ is
the set of lattice points in a polytope of the form $s \cdot {\cal P}$, 
where ${\cal P}$ is any lattice polytope and $s > {\rm dim}({\cal P})$.
\end{corollary}

\begin{proof}
Reflexive polytopes possess exactly one interior point.
If $s$ is bigger than the ambient dimension then $s$ times
any lattice polytope contains an interior point.
\end{proof}

\section{Low dimensions}

In this section we present the complete classification
of gkz-rational configurations for $d \leq 4$.
Note that the $d=1$ case is trivial since
we disallow repeated points. If we did allow them then 
$A = ( \, 1 \,\, 1 \,\,  1 \,\, \cdots \,\, 1 \,)\,$ 
would be gkz-rational for $s \geq 2$ because the function
$\,1/(x_1 + x_2 + \cdots + x_s)\,$ is $A$-hypergeometric.

Toric curves $(d=2)$ are never gkz-rational.
This was shown in  \cite[Theorem~1.10]{cdd}. We 
rederive this result as follows. Write the configuration as
$$
A\ =\ \left(
\begin{array}{cccc}
1 & 1 & \cdots & 1\\
k_1 & k_2 & \cdots & k_s
\end{array}
\right)\ ;\quad k_1<k_2<\cdots<k_s\,.
$$
Every circuit  $Z\subseteq A$ consists of three collinear points:
$$
Z\ =\ \left(
\begin{array}{ccc}
1 & 1 & 1\\
k_a & k_b &  k_c
\end{array}
\right)\ ;\quad k_a<k_b<k_c\,.
$$
Such a $1$-dimensional circuit is never balanced. Theorem~\ref{MainThm} implies
that $A$ is not gkz-rational. In what follows we 
prove the only-if part of Theorem \ref{secondmainthm}.

\begin{theorem}
\label{LowDimThm}
Let $A$ be an integer matrix with $d \leq 4$ rows.
If $A$ is gkz-rational then $A$ is affinely isomorphic
to an essential Cayley configuration.
\end{theorem}

\begin{proof}
It suffices to prove the following two assertions:
\item{$\bullet$}
If $A$ is a configuration on the line $(d=2)$ or in
$3$-space $(d=4)$ then $A$ is not gkz-rational.
\item{$\bullet$}
If $A$ is a configuration in the  plane $(d=3)$
then $A$ is gkz-rational
if and only if the points of $A$ lie on two parallel lines
with each line containing at least two points from $A$.

\vskip .1cm

The case  $d=2$ was proved above. We first assume $d=3$.
If the points of $A$ lie on two parallel lines then we can write
their coordinates as follows:
\begin{equation}
 A \quad = \quad 
\left(
\begin{array}{cccccccc}
1 & 1 &\cdots& 1 & 0 & 0 &\cdots& 0\\ 
0 & 0 &\cdots& 0 & 1 & 1 &\cdots& 1\\ 
0 & k_1 &\cdots& k_m & 0 & \ell_1 &\cdots& \ell_n\\ 
\end{array}
\right)
\end{equation}
Thus $A$ is the Cayley configuration of two one-dimensional configurations.
The construction in the next section shows
that $A$ is gkz-rational for $m,n \geq 1$.

Conversely, suppose that $A$ does not lie on 
two parallel lines. We may further assume that $A$ contains no
unbalanced spanning circuit by
Theorem \ref{MainThm}. One example of a
configuration satisfying these requirements is
\begin{equation}\label{eq:5.41}
A\  = \  \left(
\begin{array}{cccccc}
2 & 1 & 0 & 1 & 0 & 0\\ 
0 & 1 & 2 & 0 & 1 & 0\\
0 & 0 & 0 & 1 & 1 & 2
\end{array}
\right)\, .
\end{equation}
The toric variety $X_A$ is the Veronese surface in $\CP^5$.
Its dual variety is the hypersurface defined by
the discriminant of a ternary quadratic form
\begin{equation}\label{eq:5.42}
D_A \quad  = \quad   \det\left(
\begin{array}{ccc}
2x_1 & x_2 & x_4\\ 
x_2 & 2x_3 & x_5 \\
x_4 & x_5 & 2x_6 
\end{array}
\right)\,.
\end{equation}
Suppose there exists a rational $A$-hypergeometric function
$f(x) = P(x)/Q(x) $ with $Q$ a multiple of $D_A$.
Let $A'$ be the configuration obtained from $A$ by removing the
fourth and fifth columns.
Setting $x_4=x_5=0$ in $D_A$ yields 
$\,(4 x_1 x_3-x_2^2) \cdot x_6 $. We
can argue as in the proof of Theorem~\ref{MainThm} and construct
a rational $A'$-hypergeometric function whose denominator
contains the 
binomial factor. Proposition  \ref{Weakk} would imply that the 
configuration consisting of the first three columns of $A$ is gkz-rational,
and this is a contradiction to Theorem \ref{th:2.1}.
Hence the configuration $A$ in (\ref{eq:5.41}) is not gkz-rational.

Another configuration to be considered is
\begin{equation}\label{Weddge}
A\  = \  \left(
\begin{array}{cccccc}
1 & 1 & 1 & 1 & 1 \\
0 & p & q & 0 & 0 \\
0 & 0 & 0 & p & q 
\end{array}
\right),
\end{equation}
where $1 \leq p \leq q$ are relatively prime integers.
The only spanning circuit of $A$ is the balanced circuit
$\{a_2,a_3,a_4,a_5\}$. Consider the subset
$\,A' = \{a_1,a_2,a_3\}$ which is an unbalanced circuit
on the boundary of $A$. The $A$-discriminant is an irreducible
homogeneous polynomial
of degree $q^2-p^2$ which looks like 
$$ D_A \quad = \quad
x_5^{p(q-p)} \cdot D_{A'}(x_1,x_2,x_3)^{q-p} \, + \, 
\hbox{terms containing} \,\,\, x_4. $$
Applying the expansion technique with respect to $x_4$, 
we get a rational $A'$-hypergeometric function whose denominator
contains $D_{A'}$. This contradicts Theorem \ref{th:2.1}.
Hence the configuration $A$ in (\ref{Weddge}) is not gkz-rational.

Our assertion for $d=3$ now follows from the subsequent
lemma of combinatorial geometry. Note that four points in the
plane, with no three collinear, lie on two parallel lines if and
only if they form a balanced circuit.
\end{proof} 

\begin{lemma}
Let $B$ be a planar configuration without interior points such that
every four-element subset of $B$ lies on two parallel lines.
Then either $B$ lies on two parallel lines, or $B$ is affinely
equivalent to (\ref{eq:5.41}) or (\ref{Weddge}) or
to the vertices of a regular pentagon, in which case
$B$ has irrational coordinates.
\end{lemma}

\begin{proof}
We may assume without loss of generality that the origin
$O$ lies in $B$ and is a vertex of the convex hull ${\rm conv}(B)$.
Let $a$ and $b$ be the points of $B$ closest to $O$ along the edges
of ${\rm conv}(B)$ adjacent to $O$.  Let $c = a + b$.
Any other point $x\in B$ must be of the form
$ r_1\cdot a$, or $r_2\cdot b$, or $a + r_3 \cdot b$, or
$b + r_4 \cdot a$, where $r_1,r_2,r_3,r_4$ are positive real
numbers and $r_1,r_2>1$.

If $c\in B$, then only the first two cases may occur.  Indeed,
suppose
 $x= a + r_3 \cdot b\in B$, then either all the points lie on two
parallel lines or there exists a point $y = r_1\cdot a$ or
$y = b + r_4 \cdot a$ in $B$.  It is easy to check that
in all of these cases, $B$ has an interior point. Suppose then that
 $x_1 = r_1\cdot a\in B$ and
$x_2 = r_2\cdot b\in B$. We have   $r_1=r_2$ since
the subset $\{a,b,x_1,x_2\}$ lies on two parallel lines. 
If $r_1 \not= 2$ then  the subset
$\{O,c,x_1,x_2\}$ contradicts the assumption.  Hence, if 
$c\in B$, either all the points lie on two parallel lines, or 
$r_1 = r_2 = 2$ which means that $B$ is
affinely equivalent to (\ref{eq:5.41}).

On the other hand, if $c\not\in B$ and there exists a point
$x_1 = a+r_3 b\in B$, then either all the points lie on two
parallel lines or $B$ contains a point of the
form $x_2=r_1 a$ or $x_2 = b+r_4 a$. 
Since $r_3 \not= 1$, 
in all of these cases $B$ contains an unbalanced circuit,
or $B = \{a,O,b,x_1,x_2\}$ is affinely equivalent to the
vertex set of a regular pentagon.
The only remaining possibility is that all points of $B$ be
multiples of either $a$ or $b$.  But if $x_1=r_1 a $ and 
$x_2=r_2 b$ are in $B$ then $\{a,b,x_1,x_2\}$ is unbalanced
unless $r_1=r_2$.  Hence the only possible configuration not
containing the point $c$ is affinely equivalent to (\ref{Weddge})
\end{proof}

We note that the argument in the paragraph following
(\ref{eq:5.42}) works also for the
discriminant (\ref{DiscrQuadForm}) of any quadratic form.
 Hence $A' = 2 \cdot \Delta_q$ is
not gkz-rational for any $q$.
Since $A'$ is a spanning subconfiguration
of $A = r \cdot \Delta_q$ for all $r \geq 2$,
we conclude the following result which was stated 
in the introduction.

\begin{prop}
\label{MultSimp}
Multiples of simplices, $A = r \cdot \Delta_q$,
are never gkz-rational. 
\end{prop}

\vskip .1cm

We now proceed to discuss configurations in affine
$3$-space $(d=4)$. Let us begin by stating the relevant
fact of combinatorial geometry in this case.

\begin{lemma}
\label{OnlyPrismOcta}
Let $B$ be a $3$-dimensional point configuration 
which is not a pyramid and such that
every $2$-dimensional circuit is balanced
and no $3$-dimensional circuit exists.
Then either $B$ lies on three parallel lines, or 
$B$ is affinely equivalent to a subconfiguration of  
$\{ O, P,Q,R,  c P, c Q, c R \}  $
for some points $P,Q,R$ and some $c \in \R $.
\end{lemma}

\begin{proof}
Choose five points from our configuration which are not in a plane. They have 
the form $\{A_1, A_2, B_1, B_2, C\}$ 
where the lines $\overline{A_1 A_2}$ and $\overline{B_1 B_2}$ are parallel 
and $C$ lies outside the plane 
$\Pi = \overline{A_1  A_2 B_1 B_2}$. Suppose that our
configuration is not on three parallel lines. There exists a point 
$D \notin \Pi$
such that the line $\overline{C D}$ is not parallel to the lines 
$\overline{A_1 A_2}$ and $\overline{B_1 B_2}$. If, under this 
hypothesis, the line $\overline{C D}$ is still parallel to the plane
$\Pi$,  then we have created a 
$3$-dimensional circuit, a contradiction. 
Therefore the line $\overline{C D}$ meets the plane 
$\,\Pi \,$ in a point which we call the origin $O$. 
The origin $O$ must be equal to either
$\,\overline{ A_1 B_1} \cap \overline{A_2 B_2}\,$ or
$\,\overline{ A_1 B_2} \cap \overline{A_2 B_1}$; otherwise
we would have created a $3$-dimensional circuit.
{}From this requirement we conclude that the
configuration $\{O,A_1,B_1,C,A_2,B_2,D\}$ is 
affinely equivalent to
$\{ O, P,Q,R,  c P, c Q, c R \}  $.

It remains to be seen that $O$ is the only
point that may be added to the configuration
$\{ P,Q,R,  c P, c Q, c R \}  $ 
without creating either a $3$-dimensional circuit or an
unbalanced $2$-dimensional circuit.
A point not a multiple of $P,Q$ or $R$ obviously creates
a $3$-dimensional circuit. A multiple of $P,Q$ or $R$ creates
an unbalanced $2$-dimensional circuit, unless it is the
origin.
\end{proof}

\begin{proof}[Proof of Theorem \ref{LowDimThm} (continued)]
Let $A$ be a configuration in affine $3$-space. We shall
prove that $A$ is not gkz-rational. In view of Theorem \ref{MainThm},
we may assume that $A$ contains no unbalanced spanning circuit. This
implies that $A$ contains no $3$-dimensional circuit, 
because such a circuit involves five points and,
five being an odd number,  that circuit would be unbalanced.

Suppose that $A$ contains an unbalanced $2$-dimensional
circuit $Z$. Then $Z$ lies in a facet of $A$. 
There must be at least two distinct
points $P$ and $Q$ of $A$ which do not lie in that facet. 
Otherwise, $A$ is a pyramid
and the $A$-discriminant is $1$. If the line 
$\,\overline{P Q}\,$ is parallel to the plane spanned by $Z$
then,  since $Z$ is unbalanced, some triangle in $Z$ has all
of its three edges skew to $\overline{P Q}$. 
This triangle together with $P$ and $Q$
forms a $3$-dimensional circuit, a contradiction. Hence the line 
$\overline{ P Q} $ intersects the plane  spanned by $Z$.
Some triangle in $Z$ has the property that
none of the lines spanned by its edges contains that intersection point.
Again, this triangle together with $P$ and $Q$ forms a $3$-dimensional circuit.

We conclude that A has no $3$-dimensional circuit and 
every $2$-dimensional circuit of $A$  is balanced. Lemma 
\ref{OnlyPrismOcta} tells us what the possibilities are.
If $A$ lies on three parallel lines, then
$D_A = 1$ and thus $A$ is not gkz-rational.
It remains to examine the special configurations
$\,\{ O , P, cP, Q, c Q, R,  c R  \} $. An affine transformation
moves the points $P,Q$ and $R$ onto the coordinate axes, so that
our configuration has the matrix form
\begin{equation}
\label{SevenPoints}
  \left(
\begin{array}{rrrrrrr}
 1 &1 & 1 & \ 1 & 1 & \ 1 & 1  \\ 
0 & q & -p & 0 & 0 & 0 & 0    \\
0 & 0 & 0 & q & -p & 0 & 0     \\
0 & 0 & 0 & 0 & 0 & q & -p   
\end{array}
\right) 
\end{equation}
where $p$ and $q$ are relatively prime integers, and $q > 0$.
The subconfiguration consisting of the last six columns
is spanning. We shall prove that it is non-degenerate and not
gkz-rational. Our usual deletion technique then implies that
the bigger configuration (\ref{SevenPoints}) is also not gkz-rational.
It therefore suffices to consider the following $4 \times 6$-matrix
\begin{equation}
\label{SixPoints} A \quad = \quad
  \left(
\begin{array}{rrrrrrr}
1 & 1 & \ 1 & 1 & \ 1 & 1  \\ 
q & -p & 0 & 0 & 0 & 0     \\
0 & 0 & q & -p & 0 & 0     \\
0 & 0 & 0 & 0 & q & -p    
\end{array}
\right) \, .
\end{equation}
We shall distinguish the two cases $p > 0$ and $p < 0$.
If $p > 0$ then $A$ represents an octahedron, and 
if $p < 0$ then $A$ represents a triangular prism.

We shall present a detailed proof for the octahedron case $p > 0$.
The proof technique to be employed was shown to us by Laura Matusevich.
We first note that the $A$-discriminant $D_A$ is a homogeneous
irreducible polynomial of degree $(p+q)^2$.
The Newton polytope of $D_A$ is a simplex with
vertices corresponding to the monomials:
$$ (x_1^p x_2^q)^{p+q}\,,\quad (x_3^p x_4^q)^{p+q}\,,\quad
(x_5^p x_6^q)^{p+q}\, . $$
Suppose $A$ is gkz-rational.  There is a
rational function $f(x)=P(x)/Q(x)$, 
with $P,Q$ relatively prime polynomials, such that
the $A$-discriminant $D_A$ divides $Q$,
and $f$ is  $A$-hypergeometric of some degree $\beta$.  For $u \in \N^6$,
the derivative $\partial^u f$ is $A$-hypergeometric of degree
$\beta - A \cdot u$ and has $D_A$ in its denominator. 
Replacing $f$ by $\partial^u f$ for suitable $u $, we may assume
that the $A$-degree of $f$ is of the form $\beta=(K,0,0,0)$ for some
negative integer $K$.

We expand $f$ around the vertex $(x_5^p x_6^q)^{p+q}$
of the Newton polytope of $D_A$.  This results in a convergent
Taylor series in the new variables
$$
u \ :=\ (-1)^{p+q}\,\frac {x_1^p x_2^q}{x_5^p x_6^q} \quad {\hbox { and }}
\quad v \ :=\ (-1)^{p+q}\,\frac {x_3^p x_4^q}{x_5^p x_6^q}\,.
$$
That hypergeometric series equals, up to a constant,
$$ \frac {1}{(x_5)^{pk}(x_6)^{qk}} \sum_{m,n\geq 0} 
\frac {(p(m+n+k)-1))! (q(m+n+k)-1))!}
{(np)!(nq)!(mp)!(mq)!} \, u^m v^n\,,$$
for an appropriate positive integer $k$.  
The coefficients of this series can be derived directly from
the toric operators (\ref{higherorder}) arising from $A$. It is one of 
the canonical series described for general $A$ in \cite[\S 3.4]{sst2}.
The series
\begin{equation} \label{octaseries}
\Psi(u,v) \quad = \quad
\sum_{m,n\geq 0} 
\frac {(p(m+n+k)-1))! (q(m+n+k)-1))!}
{(np)!(nq)!(mp)!(mq)!} \, u^m v^n
\end{equation}
represents a rational function in two variables.

We denote by $F(m,n)$ the coefficient of $u^m v^n$
in the series  (\ref{octaseries}). Note that $F(0,0) \not= 0$.
Since $\Psi$ is rational, 
there exist positive integers
$N,M$, and constants
$c_{ij} \in \C, \; 0 \leq i,j \leq N$, such that $c_{00}\not=0$ and
$$\sum_{i,j=0}^N c_{ij} F(m+i,n+j) 
\quad =  \quad 0 \; \; \; \mbox{holds for all} \; m,n
\geq M.
$$
If we divide  $F(m+i,n+j) $ by $F(m,n)$ then we get a rational
function in $m$ and $n$. Hence the following is an identity
of rational functions in $m$ and $n$:
\begin{equation}
\label{ide}
 \sum_{i,j=0}^N c_{ij} \frac{F(m+i,n+j)}{F(m,n)} \quad = \quad 0 \; \; \;
\end{equation}

Let $R(m,n)$ and $S(m,n)$ denote  the incremental quotients:
$$R(m,n) \ :=\ \frac {F(m+1,n)}{F(m,n)}\ ;\quad
S(m,n) \ :=\ \frac {F(m,n+1)}{F(m,n)}\,.$$
If $a,b\in \N$ and we set $\mu=m+n$, $c=a+b$, we have
\begin{equation}\label{quot1}
\!\!\!\! R(m+a,n+b)\  = \ 
\frac {\prod_{j=0}^{p-1} (p(\mu+c+k)+j)  \prod_{j=0}^{q-1} 
(q(\mu+c+k)+j)}
{\prod_{j=1}^{p} (p(m+a)+j)  \prod_{j=1}^{q} 
(q(m+a)+j)} \, ,
\end{equation}
\begin{equation}\label{quot2}
\!\!\!\! S(m+a,n+b)\  = \ 
\frac {\prod_{j=0}^{p-1} (p(\mu+c+k)+j) \prod_{j=0}^{q-1} 
(q(\mu+c+k)+j)}
{\prod_{j=1}^{p} (p(n+b)+j) \prod_{j=1}^{q} 
(q(n+b)+j)}\, .
\end{equation}
Given now $0 \leq i,j \leq N$ with $i+j>0$ we have
$$
 \frac{F(m+i,n+j)}{F(m,n)}\ 
 =\  \prod_{a=0}^{i-1}  R(m+a,n+j) \cdot  \prod_{b=0}^{j-1} S(m,n+b)
$$
Note that either $R(m,n)$ or $S(m,n)$ is
a factor in the above product.  Consider now the point
$$(m_0,n_0) \ :=\ (- \frac {p-1}{p} - k - \alpha, \alpha)\,,$$
where $\alpha$ is an irrational number.  We have
$p(m_0 + n_0 + k ) = -(p-1)$ and therefore both $R(m,n)$ and $S(m,n)$
vanish at $(m_0,n_0)$.  On the other hand, since $\alpha$ is irrational,
none of the denominators in (\ref{quot1}) or (\ref{quot2}) may vanish
at $(m_0,n_0)$. Evaluating the left-hand side of (\ref{ide})
at $(m_0,n_0)$ yields $c_{00} = 0$ which is impossible.

We have shown that the matrix $A$ in (\ref{SixPoints})
is not gkz-rational for $p > 0$.
The  proof of non-rationality in the triangular prism  case $(p < 0)$,
provided to us by Laura Matusevich,  is analogous and will be omitted here.
In summary, we conclude that every 
$3$-dimensional configuration is not gkz-rational.
\end{proof}

\section{Toric residues}

In this section we present an explicit construction of non-Laurent
rational hypergeometric functions. This will prove the 
if-direction of Conjectures \ref{conj4} and \ref{conj4prime}
as promised in Theorem \ref{secondmainthm}. At the end of Section 5
we state further open problems concerning
residues and rational hypergeometric functions.
We begin with the ``Cayley trick'' for
representing resultants as discriminants.

\begin{prop}
\label{FivePointOne}
Let $A$ be a Cayley configuration (\ref{Cayley}).
If $A$ is essential then the resultant
$R_{A_0,\dots,A_r}$ is non-constant
and equals the discriminant $D_A$.
\end{prop}

\begin{proof}
The identity $\,R_{A_0,\dots,A_r} = D_A \,$
was proved in \cite[Proposition 9.1.7]{gkzbook}
under the more restrictive hypothesis that 
the configurations $A_0,\dots,A_r$ are all full-dimensional; 
see \cite[Hypothesis (1) on page 252]{gkzbook}.
The argument given in that proof shows that
$\,R_{A_0,\dots,A_r} \not= 1\,$ suffices to imply
$\,R_{A_0,\dots,A_r} = D_A $.
On the other hand, the condition of $A$ being essential
appears in \cite[equation (2.9)]{pedersen}, and 
\cite[Corollary 2.4]{pedersen} shows that it
is  equivalent to $\,R_{A_0,\dots,A_r} \not= 1$.
\end{proof}

\begin{corollary}\label{5.2}
Conjecture \ref{conj4} and 
Conjecture \ref{conj4prime} are equivalent.
\end{corollary}

\begin{proof}
We must show that a non-degenerate configuration $B$ is affinely 
isomorphic to an essential Cayley configuration
if and only if its discriminant $D_B$ equals the mixed resultant
$R_{A_0,\dots,A_r}$ of some tuple of configurations
$(A_0,\dots,A_r)$. The only-if direction is the content of
Proposition \ref{FivePointOne}. For the converse, suppose
$\,D_B = R_{A_0,\dots,A_r} \not= 1$. Let $A$ be the
Cayley configuration of $A_0,\dots,A_r$. Then
$\,D_A = D_B$. In other words, 
the toric varieties $X_A$ and $X_B$ in $\CP^{s-1}$ have the
same dual variety, namely, the hypersurface defined by
$D_A = D_B$. The Biduality Theorem \cite[Theorem 1.1.1]{gkzbook}
shows that $X_A = X_B$, and this implies that
$A$ and $B$ are affinely isomorphic.
\end{proof}

We next review the construction of the {\it toric residue}
associated with a toric variety $X_\Delta$. This was
introduced by Cox \cite{cox} and further developed
in \cite{ccd}, \cite{global}, \cite{tokyo}.
Here $\Delta$ is the set of all lattice points in a full-dimensional
convex polytope in $\R^r$.  We consider $r+1$ 
Laurent  polynomials $f_0,f_1,\dots,f_r$
supported in $\Delta$ with generic complex coefficients:
\begin{equation}\label{laurentpolys}
f_j(t) \ =\ \sum_{m\in\Delta} x_{jm}\, t^m\ ,
\qquad j=0,1,\dots,r\,.
\end{equation}
Given an interior
lattice point $a\in {\rm Int}((r+1)\cdot \Delta)$ and an index
$i\in \{0,\dots,r\}$, consider
the total sum of Grothendieck point residues:
\begin{equation}\label{globalres}
{\rm Res}^\Delta_f(t^a) \,\, := \,\,
(-1)^i\sum_{\xi\in V_i} {\rm Res}_\xi\left(
\frac {t^a/f_i}{f_0\cdots f_{i-1}f_{i+1}\cdots f_r}\,
\frac {dt_1}{t_1}\wedge\cdots\wedge \frac {dt_r}{t_r}
\right)\,,
\end{equation}
where 
$V_i \ =\ \{t\in (\C^*)^r \,:\,
f_0(t)=\cdots =f_{i-1}(t)=f_{i+1}(t)=\cdots= f_r(t)=0\}\,.$
It is shown in \cite[Theorem~0.4]{ccd} that the expression
(\ref{globalres}) is independent of $i$ and agrees with the
 residue defined by Cox in \cite{cox}. 
We refer to \cite[\S 6]{cox} and \cite[\S 5]{ccd}
for integral representations such as (\ref{anintegral})
of the toric residue $\,{\rm Res}^\Delta_f(t^a)$.

The toric residue $\,{\rm Res}^\Delta_f(t^a)\,$ is a rational function
in the coefficients $x_{jm}$ of our system (\ref{laurentpolys}).
For degree reasons, this rational function is never 
a non-zero polynomial. It was shown
in \cite[Theorem~1.4]{tokyo} that the product
$$\RR_\Delta(f_0,\dots,f_r)\cdot {\rm Res}^\Delta_f(t^a)$$
is a polynomial in the variables $x_{jm}$, where 
$\RR_\Delta(f_0,\dots,f_r)$ denotes the (unmixed) sparse resultant
associated with the polytope $\Delta$; see \cite[\S 8.2]{gkzbook}.

There is an easy algebraic method \cite[Algorithm 2]{global}
for computing  the rational function $\,{\rm Res}^\Delta_f(t^a)$:
translate $f_0,\dots,f_r,t^a$ into multihomogeneous polynomials
$F_0,\dots, F_r, u^{a'}$ in the homogeneous coordinate ring
of $X_\Delta$, compute any Gr\"obner basis
${\cal G}$ for $\langle F_0,\dots,F_r \rangle$, and finally take
the normal form modulo ${\cal G}$ of $u^{a'}$ and divide it by
the normal form modulo ${\cal G}$  of the 
{\it toric Jacobian} \cite{cox} of $F_0, \dots, F_r $.
This computation yields $\,{\rm
Res}^\Delta_f(t^a)\,$  up to a constant.

\begin{example}
We demonstrate the algorithm of \cite{global}
by showing how it computes the
rational function (\ref{anintegral}).  Here $r=1$ and $\Delta$
is the segment $[0,2]$ on the line. The  system (\ref{laurentpolys})
consists of two quadratic polynomials 
$\,  f_0(t) =  x_1  +  x_2 t +  x_3 t^2   \,$ 
and $\,  f_1(t) =  x_4  +  x_5 t +  x_6 t^2   $.
The toric variety $X_\Delta$ is the projective line
$\CP^1$. We rewrite our input equations in
homogeneous coordinates,
$$   F_0(u_1,u_2) \,= \, x_1 u_1^2 +  x_2 u_1 u_2 +  x_3 u_2^2   
\quad \hbox{and} \quad
 F_1(u_1,u_2) \,= \, x_4 u_1^2 +  x_5 u_1 u_2 +  x_6 u_2^2   , $$
and we compute any Gr\"obner basis ${\cal G}$ for the ideal
$\, \langle F_0,F_1 \rangle \,$ in ${\bf K}[u_1,u_2]$,
where $\, {\bf K}\,= \, \Q(x_1,x_2,x_3,x_4,x_5,x_6)$.
Here the toric Jacobian equals
\begin{eqnarray*}
& J(u_1,u_2) \quad = \quad
\frac{\partial F_0}{\partial u_1} \frac{\partial F_1}{\partial
u_2} \, - \,
\frac{\partial F_0}{\partial u_2} \frac{\partial F_1}{\partial
u_1} \\ &
= \quad
 2  \,(x_1 x_5-x_2 x_4) \,u_1^2
\,\, + \,\,4 \,(x_1 x_6- x_3 x_4) \,u_1 u_2
\,\, + \,\, 2 \, ( x_2 x_6- x_3 x_5) \,u_2^2.
\end{eqnarray*}
The residue $\,  {\rm Res}^\Delta_f( t^2 ) $, which appears in
(\ref{anintegral}), is computed
as  4 times the ratio of the normal form of $u_1 u_2$
over that  of $J(u_1,u_2)$, both modulo ${\cal G}$.
\end{example}

To establish the connection to hypergeometric functions, we
now consider the Cayley configuration of $\Delta,\Delta,\dots,
\Delta$, taken $r+1$ times:
$$ \widehat\Delta \quad := \quad
\bigcup_{i=0}^r ( \{e_i\}\times \Delta )
\,\,\, \subset \,\,\,
\Z^{r+1} \times \Z^{r} \,\,= \,\, \Z^{2r+1}\,. $$
The points in $\widehat \Delta$ are labeled
by the variables $x_{im} $ for $i=0,\dots,r $ and $m \in \Delta$.

\begin{lemma}
\label{AllSame}
The  configuration $\widehat\Delta$ is gkz-rational.
\end{lemma}

\begin{proof}
Let $\,a \in {\rm Int}((r+1)\cdot \Delta)$ and $f_0,\dots, f_r$ generic
Laurent polynomials as in (\ref{laurentpolys}).
It follows from either the definition or
\cite[Algorithm 2]{global} that ${\rm Res}^\Delta_f(t^a)$ 
is a homogeneous function with respect to the grading 
induced by $\widehat\Delta$; that is, it satisfies 
the equations defined by
the operators (\ref{firstorder}) for a suitable parameter vector.
It follows from \cite[Theorem~7]{global} that
${\rm Res}^\Delta_f(t^a)$ is also annihilated by the
 operators (\ref{higherorder}).
Hence $\,{\rm Res}^\Delta_f(t^a)\,$ is a rational
$\widehat \Delta$-hypergeometric function.

The discriminant associated with $\widehat \Delta$ equals the resultant
$\,\RR_\Delta =
R_{
\Delta,
\Delta,\dots,
\Delta} $, by Proposition \ref{FivePointOne}.
This resultant is not a monomial, for instance, by
\cite[Corollary 8.2.3]{gkzbook}. We showed in 
\cite[Theorem~1.4]{tokyo} that
$\, \RR_\Delta  \cdot {\rm Res}^\Delta_f(t^a)\,$
is a polynomial. It remains to be seen that the
toric residue $\, {\rm Res}^\Delta_f(t^a)\,$ itself
is non-zero
for at least one lattice point $\,a \in {\rm Int}((r+1)\cdot \Delta)$.
Recall from \cite[Theorem 5.1]{cox} and \cite[Proposition 1.2]{tokyo}
that the polynomial 
$$
j(t) \ =\ \det \left( 
\begin{array}{cccc}
f_0 &
f_1 & \dots &
f_r  \\
t_1\frac {\partial f_0}{\partial t_1} &
t_1 \frac {\partial f_1}{\partial t_1} & \dots &
t_1 \frac {\partial f_r }{\partial t_1}  \\
\vdots & \vdots & \ddots & \vdots \\
t_r \frac {\partial f_0}{\partial t_r} &
t_r \frac {\partial f_1}{\partial t_r} & \dots &
t_r \frac {\partial f_r }{\partial t_r}  \\
\end{array}
\right)
$$
is supported in ${\rm Int}((r+1)\cdot \Delta)$
and ${\rm Res}^\Delta_f(j(t))  =  n!\cdot{\rm vol}(\Delta) $.
Here the operator
$\,{\rm Res}^\Delta_f( \,\cdot\,) \,$ is extended
from monomials $t^a$ to the polynomial $j(t)$ by
linearity. At least one
of the residues $\,{\rm Res}^\Delta_f(t^a) \,$
 as $a$ runs over $\,{\rm Int}((r+1)\cdot \Delta)$,
does not vanish and hence is a non-Laurent
rational $\widehat\Delta$-hypergeometric function.
\end{proof}

\begin{example}
\label{determinant}
The reciprocal of the determinant is a hypergeometric function.
To see this, take $\Delta$ to be the unit simplex, so that,
$\,f_i =  x_{i0} + 
x_{i1} t_1 +  \cdots + x_{ir} t_r \,$
in affine coordinates on $X_\Delta = \CP^r$.
The scaled simplex $\,(r+1)\cdot\Delta \,$ has a unique interior
lattice point $a$, and
$\, {\rm Res}^\Delta_f(t^a) \,=\, 1 / {\rm det}(x_{ij})$.
Here the Cayley configuration $\widehat \Delta $ 
is the product of two $r$-simplices $ \Delta \times \Delta $.
We conclude that  $\, 1 / {\rm det}(x_{ij})\,$ is
a rational  $ \Delta \times \Delta $-hypergeometric function.
\end{example}

\vskip .1cm

We are now prepared to complete the proof of our main result.

\begin{proof}[Proof of the if direction in Theorem \ref{secondmainthm}]
Let $A$ be an essential Cayley configuration with 
$A_0,A_1,\dots,A_r$ as in (\ref{Cayley}). Let $\Delta$
 be the set of all lattice points in a convex polytope
containing all the configurations 
$\,A_i $ for $i=0,1,\dots,r$. Then, $\Delta$ is full-dimensional and
$\,A \subseteq \widehat \Delta$.

Consider  configurations $B_0,\dots,B_r$ in  $\Z^r$ such
that $A_i\subseteq B_i \subseteq \Delta$ for $i=0,\dots,r$.
The corresponding  Cayley configuration $B$ is
still essential, since the Minkowski sum
$\,\sum_{i\in I} B_i\,$ has affine dimension
at least $\,|I|$. This property holds for $A$ 
and it does for $B$. We conclude from
Proposition  \ref{FivePointOne} that 
the Cayley configuration $B$ is
non-degenerate and $D_B = R_{B_0,\dots,B_r}$.

We would like to show that, in fact, any such configuration
$B$ must be gkz-rational.  We proceed by induction on the
cardinality of $\widehat \Delta \backslash B$. The base case
is cardinality zero:  if $B = \widehat \Delta $
then $B$ is gkz-rational by  Lemma \ref{AllSame}.

For the induction step we may  suppose that
$\tilde B$ is obtained
from $B$ by removing a point $b$ from $B_0 \backslash A_0$ and
assume, inductively, that $f$ is a rational $B$-hypergeometric
function which contains the discriminant $D_B$ in its denominator.
Let us denote by 
$t$ the variable associated with $b$ and by $\tilde x$ the variables
associated with $\tilde B$. 
Expand as
in (\ref{SeriesForMircea}):
\begin{equation}
\label{SeriesForMircea2}
f(\tilde x;t) \quad = \quad 
\sum_{\ell\geq \ell_0}\,R_\ell(\tilde x)\cdot t^\ell\,,
\end{equation}
where each $R_\ell(\tilde x)$ is a rational $\tilde B$-hypergeometric
function.  
We may now argue as in the proof of Theorem~\ref{MainThm}; since
$B$ and $\tilde B$ have affine dimension $2r$ it follows from
Lemma~\ref{setzero} that the $\tilde B$ discriminant $D_{\tilde B}$
divides the specialization $D_B|_{t=0}$.  Hence, for some $\ell$,
the rational function $R_\ell(\tilde x)$ will lie strictly in the
field of fractions of the domain $\C[\tilde x,\tilde x^{-1}]_{
\langle D_{\tilde B} \rangle}$ and, consequently, will be a rational
$\tilde B$-hypergeometric function which contains  the discriminant
$ D_{\tilde B}$ in its denominator.
In summary, the configuration $\, \tilde B  \,= \,
B \backslash \{b\}$ inherits the property of being
gkz-rational from the configuration $B$. By induction,
we conclude that $A$ is gkz-rational.
\end{proof}

The results in this paper raise many questions about
rational hypergeometric functions.
The most obvious one is whether 
Conjectures \ref{conj4} and \ref{conj4prime} are true for
toric varieties other than hypersurfaces, curves, surfaces and threefolds.
Another question which concerns the number of rational solutions is the
following:
Is the dimension of the vector space of rational function solutions
to the hypergeometric system  $H_A(\beta)$
always bounded by the normalized volume of $A$?
This volume is the degree of $X_A$; see (\ref{rankequation}).

It would be nice to extend the observation in Example \ref{determinant}
from determinants to hyperdeterminants. Following \cite[Chapter 14]{gkzbook},
the {\it hyperdeterminant} is the discriminant (= dual hypersurface)
associated with any Segre variety
$\, X_A =  \CP^{k_0} \times \CP^{k_1} \times \cdots \times \CP^{k_r} $.
Suppose $\,k_0 \geq k_1 \geq \cdots \geq k_r$. The corresponding
configuration is a product of simplices
$\,A =  \Delta_{k_0} \times \Delta_{k_1} \times \cdots \times \Delta_{k_r} $.
It is known \cite[Theorem 14.1.3]{gkzbook} that
$A$ is non-degenerate if and only if 
$\, k_0 \leq  k_1 + k_2 + \cdots + k_r $.
The case of equality 
$\, k_0 =  k_1 + k_2 + \cdots + k_r \,$ is of special interest;
it defines the {\it hyperdeterminants of boundary format}.
Since in this case $A$ is an essential Cayley configuration,
it follows from Theorem~\ref{secondmainthm} that
$A$ is gkz-rational.
We conjecture  the converse of this statement:
\begin{conjecture}
Let $A$ be the
product of simplices
$\, \Delta_{k_0} \times \Delta_{k_1} \times \cdots \times \Delta_{k_r} \,$
where $\,k_0 \geq \cdots \geq k_r $.
Then $A$ 
is gkz-rational if and only if
$\, k_0 =  k_1 + \cdots + k_r $.
\end{conjecture}

Finally, we are hoping for a 
``Universality Theorem for Toric Residues'' to the
effect that the space of rational hypergeometric functions
is spanned by Laurent polynomials and toric residues.
This statement is literally false, as the following example shows.
Let $A$ be the Cayley configuration of the segments $\{0,1\}$ and
$\{0,2\}$. The following rational function is $A$-hypergeometric:
$$ f(x_1,x_2,x_3,x_4) \quad = \quad  
\frac{x_4\,(- x_1^4x_4^2 - 6 x_1^2x_2^2x_3x_4
 +3x_2^4x_3^2)}{x_2^2 \,(x_2^2x_3+x_1^2x_4)^3}. $$
This is not a toric residue because the
degree is zero in the variables $\{x_3,x_4\}$.
However, an appropriate derivative of $f$ 
will be a toric residue. For example,
$$  \frac{\partial f}{\partial x_4} \quad  = \quad
 3 \ \frac{x_3\,(x_1^4x_4^2 - 6 x_1^2x_2^2x_3x_4
+x_2^4x_3^2)}{(x_2^2x_3+x_1^2x_4)^4} $$
agrees, up to constant, with the toric residue  in $\CP^1$
associated with the differential form
$$ \frac{t^4}{(x_1 + x_2 t)^4\cdot (x_3 + x_4 t^2)}\,
{\frac{dt}{t}}. $$

Although this is not explicit in the proof of
Theorem~\ref{secondmainthm}, one can show that
every essential Cayley configuration
admits rational hypergeometric functions 
which are toric residues and whose denominators
are multiples of the $A$-discriminant. 
Families of examples together with
extensive computer experiments support the following
general conjecture:

\begin{conjecture} \label{conj5}
Every rational  $A$-hypergeometric function $f$
has an iterated derivative
$$\frac{\partial^{i_1 + \cdots + i_s}}{
\partial x_1^{i_1} \cdots \partial x_s^{i_s}} f $$
 which is a
toric residue defined by  some facial subset  of $A$.
\end{conjecture} 

\bigskip
\bigskip

\noindent{\bf Acknowledgements:}  This work started while
the authors were visiting the Mathematical Sciences Research
Institute (MSRI Berkeley)
during the special semester (Fall 1998)
 on ``Symbolic Computation
in Geometry and Analysis". We wish to thank MSRI for its
hospitality and for the wonderful working atmosphere it provided.
We are very grateful to Mircea Musta\c{t}\v{a} for the
proof of Lemma~\ref{mircea}
and to Laura
Matusevich who devised the technique to dispose of some of the
exceptional
configurations in \S 4.  Alicia Dickenstein was partially
supported by UBACYT TX94 and CONICET, Argentina.
Bernd Sturmfels was partially supported
by NSF Grant DMS-9970254.

\bigskip
\bigskip
\medskip

{\small
\begin{flushright}
Eduardo Cattani \\
{\tt cattani@math.umass.edu} \\
Dept.~of Mathematics and Statistics \\
University of Massachusetts \\
Amherst, MA 01003, USA 
\end{flushright}

\medskip

\begin{flushright}
Alicia Dickenstein \\
{\tt  alidick@dm.uba.ar} \\
Dto.~de Matem\'atica, Universidad de Buenos Aires \\
(1428) Buenos Aires, \ Argentina 
\end{flushright}

\medskip

\begin{flushright}
Bernd Sturmfels \\
{\tt bernd@math.berkeley.edu} \\
Dept.~of Mathematics, University of California \\
Berkeley, CA 94720, USA
\end{flushright}}
\end{document}